\documentclass[11pt]{article}
\usepackage{amsmath,amssymb,amsthm,amscd,tabularx,array,
fancyhdr}
   
\newfont{\extra}{msbm10 scaled\magstep1}

\newcommand{\extr}[1]{\mbox{\extra #1}}

\newcommand{\sect}[1]{\setcounter{equation}{0}\section{#1}}
\newcommand{\subsect}[1]{\subsection{#1}}

\newtheorem{theorem}{Theorem}[section]

\def\be{\begin{equation}}
\def\ee{\end{equation}}
\def\bea{\begin{eqnarray}}
\def\eea{\end{eqnarray}}

\def\C{\extr C}
\def\K{\extr K}
\def\N{\extr N}
\def\R{\extr R}

\def\k{\kappa}

\newcommand{\bicross}{\triangleright\!\!\!\blacktriangleleft}

\newcommand{\leco}{>\!\!\blacktriangleleft}

\newcommand{\RL}{\triangleright\!\!\!\blacktriangleleft}
\newcommand{\LR}{\blacktriangleright\!\!\!\triangleleft}

\newcommand{\RIMO}{\triangleright\!\!\!<}

\newcommand{\LEMO}{>\!\!\!\triangleleft}
\newcommand{\RICO}{\blacktriangleright\!\!\!<}

\newcommand{\lact}{\triangleright}
\newcommand{\ract}{\triangleleft}
\newcommand{\lcact}{\blacktriangleleft}
\newcommand{\rcact}{\blacktriangleright}

\begin{document}

\begin{center}{ \LARGE \bf
Induced Representations \\[0.4cm]
 of Quantum   Kinematical  Algebras \\[0.4cm] and Quantum Mechanics}
\end{center}
\vskip0.25cm

\begin{center}
Oscar Arratia $^1$ and Mariano A. del Olmo $^2$
\vskip0.25cm
{ \it $^{1}$ Departamento de  Matem\'atica Aplicada a la  Ingenier\'{\i}a,  \\[0.15cm]
$^{2}$ Departamento de  F\'{\i}sica Te\'orica,\\[0.15cm]
 Universidad de  Valladolid, 
 E-47011, Valladolid,  Spain}
\vskip0.15cm

E. mail: oscarr@wmatem.eis.uva.es, olmo@fta.uva.es
 
\end{center}
 
\vskip1.5cm
\centerline{\today}
\vskip1.5cm

\begin{abstract}
Unitary representations of kinematical symmetry groups of quantum systems are fundamental in
quantum theory. We propose in this paper  its generalization to quantum kinematical groups.
Using the method, proposed by us in a recent paper
\cite{olmo01}, to induce representations of quantum bicrossproduct algebras we construct the
representations of the family of standard quantum inhomogeneous algebras
$U_\lambda(iso_{\omega}(2))$. This family contains  the quantum Euclidean, Galilei and  
Poincar\'e algebras, all of them in (1+1) dimensions. As byproducts we obtain the actions of
these quantum algebras on regular co-spaces that are an algebraic generalization of the
homogeneous spaces and $q$--Casimir equations which play the role of $q$--Schr\"odinger
equations. 
\end{abstract}
\newpage
\sect{Introduction}

It is well known  the role played by the unitary 
representations of the (kinematical) symmetry groups of quantum physical systems. One can
classify the elementary systems \cite{wigner39} according to them or to obtain their
corresponding Schr\"odinger  equations using local representations
\cite{olmo84}. In the last years we are involved in a program for constructing a theory of
induced representations for quantum algebras \cite{olmo01,olmo98,ao00,olmo00}.
  Our aim is to obtain a quantum counterpart of the  program that
Wigner started in 1939 for Lie groups  \cite{wigner39} that has been so fruitful in quantum
physics.

Kinematical groups like Poincar\'e and Galilei, whose
physical interest is out of any doubt, have a structure of semidirect product.   The quantum
version of this kind of structure for quantum groups is that of bicrossproduct 
\cite{Maj88,Maj95a}. In this case a quantum Lie algebras inherits  the `semidirect' structure
in the algebra sector and the algebra of functions also has  a semidirect product structure in
the coalgebra sector.
The interest of the quantum versions  of the kinematical groups and algebras from a physical
point of view is as $q$--generalizations of the  symmetries of the physical
space-time in a noncommutative framework. The study of  these quantum symmetries and their
representations generalizes the  Wigner program inside the perspective of the  noncommutative
geometry \cite{connes}, whose importance in  physics is increasing.

 For these groups having structure of semidirect product 
Mackey's method \cite{mackey} provides their unitary representations by induction from the
representations of one of their subgroups. In this paper we study the induced representations
for some quantum kinematical algebras from a unified point of view.  We consider the framework
of the Cayley--Klein (CK) pseudo-ortoghonal algebras \cite{olmo93} acting in spaces of
$(1+1)$--dimensions, which includes the inhomogeneous algebras like   Galilei,  Poincar\'e and 
Euclidean algebras. The quantum deformation of these CK algebras was performed in \cite{olmo93a}
(Ref.~\cite{olmo94}--\cite{olmo95} present quantum CK algebras in higher dimensions). The above
mentioned quantum CK inhomogeneous algebras have a bicrossproduct structure \cite{olmo97} that
we profit to obtain their induced representations, given a unified model for all of them.   

The induction procedure formulated by  us presents a strong algebraic character
because  we made use  of objects like modules, comodules, etc.,  which, from our point of view,
are the appropriate tools to work with the algebraic structures  characteristic of the quantum
algebras and groups. 

In the literature we can find some attempts to develop techniques that generalize the Mackey
method  \cite{mackey} of induced representations for semidirect product groups to the quantum
case. For instance, Dobrev presented in
\cite{dobrev1} a method for constructing representations of quantum groups  near
to ours, since both methods  emphasize the dual case, closer to the classical one, and the
representations are constructed in the algebra sector.  
Other authors \cite{ibort}--\cite{bgst98} have also extended 
the induction technique to quantum groups but constructing corepresentations, i.e.
representations of the coalgebra sector.   

The organization of the paper is the following one. In Section~\ref{preliminaries} we present a
mathematical  outline of the  concepts that we will use along the paper in order to unify
notation.  Section~\ref{inducedrepresentations} is devoted to summarize  the theory of induced
representations of quantum bicrossproduct  algebras. We begin to study the induction problem 
taking into account the deep relation between modules and representations obtaining, in some
sense, more deep results from a geometric point of view using the concept of regular co-space.
Since our aim is to construct the induced representations of certain quantum inhomogeneous CK 
algebras  Section~\ref{cayley-kleinalgebras} is devoted to describe briefly these $q$--algebras.
In Section~\ref{representationsquantumstandardulambda} we construct with some detail the
induced representations following the method developed
in Section~\ref{inducedrepresentations}. We start with the computation of the flow
associated to the accion of the generator of one of the  factors of the bicrossproduct over the
other one. After, we are able to determine the regular co-spaces determining the induced
representations, that we obtain below in an explicit way. Moreover, $q$--Casimir equations are
also obtained.

\sect{Preliminaries}\label{preliminaries}

Let $H=\left({V};{m}\ {\eta};{\Delta}\
{\epsilon};{S}\right)$ be a Hopf algebra with underlying vector space $V$ over the field
$\K$ ($\C$ or $\R$), multiplication 
$m: H \otimes H \to H$, coproduct $\Delta : H \to H\otimes H$, unit $ \eta:
\K \to H $, counit $ \epsilon: H \to \K$ and  antipode $S: H \to
H$. 

The algebras involved in this  work are infinite dimensional algebras but finitely
generated, for this reason  we will use a multi-index notation \cite{olmo00}.  Let $A$ be an
algebra generated by the  elements $(a_1, a_2, \ldots, a_r)$ such that the ordered monomials 
$a^n = a_1^{n_1} a_2^{n_2} \cdots a_r^{n_r}$  $( 
n=(n_1, n_2, \ldots, n_r) \in \mathbb{N}^r)
$ form a basis of the linear space  underlying to $A$.
An arbitrary product of generators of $A$ is written in a normal ordering if it is 
expressed in terms of the (ordered) basis  $(a^n)_{n \in \mathbb{N}^r}$.  The unit of $A$,
$1_A$, is denoted by  $a^0 $ ($0 \in \N^n$). Multi-factorials
and multi-deltas are defined by
 \begin{equation}
  l! = \prod_{i=1}^n l_i!, \qquad \delta_l^m= \prod_{i=1}^n \delta_{l_i}^{m_i}.
\end{equation}

A pairing between two Hopf algebras, $H$ and $H'$, is a bilinear mapping
$\langle \, \cdot \, ,  \, \cdot \, \rangle:
H \times H' \rightarrow \K$ verifying some defining relations \cite{ChP}.
The pairing is said to be left (right) nondegenerate  if
$[\langle h, \varphi \rangle = 0, \  \forall \varphi \in H'] \Rightarrow h=0$
($[\langle h, \varphi \rangle = 0, \  \forall h \in H] \Rightarrow \varphi=0$). The  pairing is
said nondegenerate if it is simultaneously  left and right nondegenerate. 
Two Hopf algebras and a  nondegenerate  pairing, $(H, H', \langle\, \cdot \, , \, \cdot \,
\rangle)$, determine a `nondegenerate triplet'. 
The bases $(h^m)$ of $H$ and
$(\varphi_n)$ of $H'$ are dual with respect to the nondegenerate pairing if
\begin{equation}
\langle h^m, \varphi_n \rangle = c_n \delta^{m}_{n},\qquad  c_n \in \K -\{ 0 \}.
\end{equation}
Given a nondegenerate triplet  and the map 
 $f: H \rightarrow H$. We shall say that     the map  $f^\dagger: {H'} \rightarrow { H'}$, 
defined  by
\begin{equation}
\langle h, f^\dagger(\varphi) \rangle= \langle f(h), \varphi \rangle ,
\end{equation}
is  the adjoint map to $f$ with respect to $\langle\, \cdot \, , \, \cdot \, \rangle$.

We will use the endomorphisms of `multiplication' (denoted
again by $h$ and $\varphi$) and `formal derivation'
($\frac{\partial}{\partial h}\equiv \partial _h$ and
$\frac{\partial}{\partial \varphi}\equiv \partial _\varphi$).
They verify that 
 $
 h^\dagger= {\partial_\varphi}$ and
 $\varphi^\dagger= {\partial_h}$.
The situation is similar when $H$ and $H'$ are finitely generated.
The definition of the formal derivative  is 
\begin{equation}
\displaystyle {\frac{\partial}{\partial h_i}} (h_1^{l_1} \cdots h_i^{l_i} \cdots
h_n^{l_n}) =
 l_i \; h_1^{l_1} \cdots h_i^{l_i-1} \cdots h_n^{l_n}. 
\end{equation}
The generalization of the `multiplication' operators is not straightforward if
the algebras are non commutative. To avoid any confusion the formal
operators associated to the generators  $h_i$ 
 will be denoted using a bar over the corresponding symbol. The
action of these operators is 
\begin{equation}
  \overline{h}_i \; (h_1^{l_1} \cdots h_i^{l_i} \cdots h_n^{l_n})=
  h_1^{l_1} \cdots h_i^{l_i+1} \cdots h_n^{l_n}. 
\end{equation}
For the elements $\varphi^i \in H'$ the `multiplication' and derivative operators are
defined in a similar way. Note that if
$H$ (resp.
$H'$) is commutative then
$\overline{h}_i$ (resp. $\overline{\varphi}^i$) acts as
a multiplication operator, but this is not longer true when
the algebra is noncommutative. 
The adjoint operators are
$\overline{h}_i^\dagger=  \partial_{\varphi^i}, \
\overline{\varphi}^{i\dagger} = \partial_{h_i}$.
The commutation relations for $\overline{h}_i$ and
$\partial_{h_i}$ (similar for $\overline{\varphi}^i$ and $\partial_{\varphi^i}$) are  
 \begin{equation}
 [\partial_{h_i}, \overline{h}_j]= \delta_{ij}, \qquad
[\overline{h}_i, \overline{h}_j]=0, \qquad
  [\partial_{h_i}, \partial_{h_j}]= 0 .
 \end{equation}


Let $(V, \alpha, A)$ be a triad composed by a unital and associative  $\mathbb{K}$--algebra
$A$, a  $\mathbb{K}$--vector space $V$ and  a
linear map (called  action) $\alpha: A\otimes_{\mathbb{K}} V \rightarrow V$ 
($\alpha(a\otimes v)=a \lact v$). It is said that  $(V, \alpha, A)$ 
(or $(V, \lact, A)$) is a left $A$--module if:
\begin{equation}
    a \lact (b\lact v)= (ab) \lact v, \qquad 1 \lact v= v, \qquad
         \forall a, b \in A, \quad \forall v \in V.
\end{equation}
Dualizing an $A$--module a comodule  is obtained.
The triad  $(V, \lcact, C)$, where  $C$ is an associative
$\mathbb{K}$--coalgebra  with counit, $V$ a  $\mathbb{K}$--vector space and 
$\lcact : V \rightarrow C \otimes_{\mathbb{K}} V$ ($v\lcact= v^{(1)} \otimes
v^{(2)}$) a (linear map) coaction,  is said to be a 
 left $C$--comodule if 
  \begin{equation} \label{comodulo}
 {v^{(1)}}_{(1)} \otimes {v^{(1)}}_{(2)} \otimes {v^{(2)}}= 
           {v^{(1)}} \otimes {v^{(2)}}_{(1)} \otimes {v^{(2)}}_{(2)}, \quad
           \epsilon(v^{(1)}) v^{(2)}= v, \qquad \forall v \in V,
\end{equation}
with  $\Delta(c)=c_{(1)} \otimes c_{(2)}$ denoting the coproduct of the elements of $C$.

A morphism of left  $A$--modules, $(V, \lact,A)$ and $(V', \lact',A)$, is a linear
map, $f:V\rightarrow V'$,  equivariant with respect the action, i.e., 
\begin{equation} 
f(a\lact v)= a \lact' f(v), \qquad \forall a \in A,\ \forall v \in V. 
\end{equation}
A linear map $f:V\rightarrow V'$ between two  $C$--comodules, $(V, \lcact ,C)$ and  
$(V',\lcact',C)$ is a morphism if
\begin{equation} \label{comodulomorfismo}
   v^{(1)} \otimes f(v^{(2)})=f(v)^{(1)'}\otimes  f(v)^{(2)'},
   \qquad  \forall v \in V.
\end{equation}

If a bialgebra acts or coacts on a vector space equipped with an additional 
structure (algebra, coalgebra or bialgebra) some  compatibility 
relations for the action may be demanded \cite{Maj95a}.   

Let $A,\ B,\ C$ be an algebra, a bialgebra and a coalgebra, respectively.
The left module  $(A, \lact, B)$ is a $B$--module algebra if
$m_A$ and $\eta_A$ are morphisms of $B$--modules. That is, if 
\begin{equation} 
b\lact( aa')= (b_{(1)}\lact a)(b_{(2)}\lact a'), \qquad 
b \lact 1 = \epsilon(b) 1, \qquad 
\forall b \in B, \, \forall a, a' \in A.
\end{equation}
A left $B$--module  $(C, \lact, B)$ is a $B$--module coalgebra if
$\Delta_C$ and $\epsilon_C$ are morphisms of $B$--modules, i.e., if 
$$
 (b\lact c)_{(1)} \otimes (b \lact c)_{(2)} =
    (b_{(1)}\lact c_{(1)}) \otimes (b_{(2)}\lact c_{(2)}), 
  \quad \epsilon_{C}(b \lact c) = \epsilon_{B}(b) \epsilon_{C}(c),
\qquad \forall b,c \in B.
$$
A left $B$--comodule  $(C, \lcact, B)$ is said to be a  $B$--comodule coalgebra if
$\Delta_C$ and $\epsilon_C$ are morphisms of $B$--comodules, i.e., 
$$
    c^{(1)} \otimes {c^{(2)}}_{(1)} \otimes {c^{(2)}}_{(2)}=
   {c_{(1)}}^{(1)} {c_{(2)}}^{(1)} \otimes {c_{(1)}}^{(2)}
    \otimes {c_{(2)}}^{(2)}, \quad
     c^{(1)} \epsilon_C (c^{(2)})= (\eta_{B} \circ \epsilon_C)(c).
$$
A left $B$--comodule  $(A, \lcact, B)$ is a $B$--comodule algebra if
$m_A$ and $\eta_A$ are morphisms of $B$--comodules. Explicitly 
\begin{equation}
(aa')_{(1)} \otimes (aa')_{(2)} = a_{(1)}a'_{(1)} \otimes  a_{(2)} a'_{(2)}, \quad  
1_A \lcact = 1_B \otimes 1_A.
\end{equation}

The triad $(B', \lact, B)$ is a left $B$--module  bialgebra if simultaneously is a
$B$--module algebra and a $B$--module coalgebra; $(B', \lcact, B)$ is a left
$B$--comodule bialgebra if  is a $B$--comodule algebra and a
$B$--comodule coalgebra.

A regular module (comodule) is an $A$--module ($C$--comodule) whose
vector space is the underlying vector space of the algebra $A$ (coalgebra $C$). The
action (coaction) is defined in terms of the algebra product (coalgebra coproduct).  
If $B$ is a bialgebra,  the regular $B$--module $(B, \lact, B)$, whose regular action is
$ b \lact b' =b b'$,
 is a module coalgebra. The regular module  $(B^*, \ract, B)$, obtained by dualization,
 is a module algebra with  regular action   
$\varphi \ract b = \langle \varphi_{(1)},b \rangle \varphi_{(2)}$ 
with $b  \in B, \; \varphi \in B^*$. 


Let $K$ and $L$ be two Hopf algebras, with  
$(L,\ract,K)$   a right  $K$--module algebra and $(K,\lcact,L)$ a left
$L$--comodule coalgebra. The tensor product
$K\otimes L$ is equipped simultaneously with the semidirect  structures of  algebra
$K\RIMO L$ and coalgebra $K \leco L$. If some  compatibility conditions
are verified  $K \RIMO L$ and  $K \leco L$ determine a Hopf algebra called (right--left) 
bicrossproduct and denoted by  $K \RL L$ \cite{Maj88,Maj95a}.

Let $\cal K$ and $\cal L$ be two  Hopf algebras and  
$({\cal L}, \lact,{\cal K})$ and $({\cal K}, \rcact,{\cal L})$
a left $K$--module algebra and a
 right $L$--comodule coalgebra, respectively, verifying
the corresponding compatibility  conditions.
Then ${\cal L} \LEMO {\cal K}$ and  ${\cal L} \RICO {\cal K}$ determine
a Hopf algebra, ${\cal L} \LR {\cal K}$, called (left--right) bicrossproduct. 

These two  bicrossproduct structures are
related by duality: if 
$K$ and $L$ are  finite dimensional bialgebras  and
the  $K$--module algebra $(L, \ract, K)$ and the 
$L$--comodule coalgebra $(K,\lcact,L)$  verify the compatibility  
conditions determining the bicrossproduct $K \RL L$ , then  $(K\RL L)^*= K^* \LR L^*$.


It has been proved in Ref.~\cite{olmo01} that  dual bases and
$*$--structures over bicrossproduct Hopf algebras may be  constructed when the 
corresponding bases and  $*$--structures of the bicrossproduct the  factors are known.
Thus, let $H=K \RL L$ be a  bicrossproduct Hopf algebra and $(K,K^*,\langle \cdot, \cdot
\rangle_1)$ and $(L, L^*,\langle \cdot, \cdot \rangle_2)$ nondegenerate triplets, 
 then the expression
 \begin{equation} \label{pairingprod}
     \langle kl, \kappa \lambda \rangle=
    \langle k, \kappa \rangle_{1}
\langle l, \lambda \rangle_{2}.
 \end{equation}
defines a nondegenerate pairing between $H$ and $H^*$. If   $(k_m)$ and
$(\kappa_m)$ are dual bases for $K$ and $K^*$, and $(l_n)$ and $(\lambda_n)$  for $L$ and
$L^*$, then
$(k_m l_n)$ and $(\kappa^m \lambda^n)$ are dual bases for $H$ and $H^*$. 

It can be proved that given a  bicrossproduct Hopf algebra, $H=K\RL L$, such that
$K$ and $L$ are equipped with $*$--structures verifying the compatibility relation
$(l\ract k)^*= l^* \ract S(k)^*$, there is a $*$--structure on the algebra sector of $H$ is
determined by
 \begin{equation} \label{estrellabicross}
     (kl)^*= l^* k^*, \qquad k \in K, \ l \in L .
  \end{equation}

\sect{Induced representations of  bicrossproduct algebras}
\label{inducedrepresentations}

The main results about the theory of induced representations for quantum bicrossproduct 
algebras will be summarized in this Section (see Refs.~\cite{olmo01,olmo98,ao00,olmo00} for a
complete description of induced representations of quantum algebras).

As we have seen in the previous Section different actions are involved, henceforth we 
will denote them by the following symbols (or their symmetric for the corresponding
right actions and  coactions):
$\lact$  ($\lcact$) (bicrossproduct actions  (coactions)),
$\vdash$ (induced and inducting representations) and
$\succ$, $\prec$  (regular actions).

Let $(H, {\cal H}, \langle \cdot, \cdot\rangle)$ be a nondegenerate triplet and
$L$ a commutative subalgebra of $H$. Suppose that  $\{ l_1, \ldots,l_s\}$ is a system
of generators of  $L$ which is completed with
$\{ k_1,\ldots,k_r\}$ to get a system of generators of $H$, such that  
$(l_n)_{n \in \mathbb{N}^s}$ is a basis of
$L$ and $(k_ml_n)_{(m,n) \in \mathbb{N}^r \times \mathbb{N}^s}$ a basis of 
$H$. Moreover,  there is a generator system of $\cal H$, 
$\{ \kappa_1, \ldots, \kappa_r, \lambda_1, \ldots,\lambda_s\}$, such that 
$(\kappa^m \lambda^n)_{(m,n) \in \mathbb{N}^r \times \mathbb{N}^s}$ is a basis
of $\cal H$ dual of that of  $H$ with pairing
\begin{equation} \label{pairing}
\langle k_m l_n, \kappa^{m'} \lambda^{n'} \rangle = m! n! \;
\delta^{m'}_m \delta^{n'}_n .
\end{equation}

In order to construct the  representation of $H$ induced by  the character
of $L$, determined by  $a=(a_1,\ldots, a_s) \in \mathbb{K}^s$ 
(explicitly, $1 \dashv l_n = a_n=a_1^{n_1} \cdots a_s^{n_s}, \ n \in  \mathbb{N}^s$) we need to
know   its carrier space  $\mathbb{K}^\uparrow$ and the action of  $H$ on it.
The  elements of $\mathbb{K}^\uparrow$ are those of
$\text{Hom}_\mathbb{K}(H,\mathbb{K})$ verifying the invariance  condition
\begin{equation}\label{equivariancecondition}
   f(hl)= f (h) \dashv l, \qquad
    \forall l \in L, \quad \forall h \in H.
\end{equation}
They can be written as
\begin{equation}
f= \sum_{(m,n)\in \mathbb{N}^r \times \mathbb{N}^s} f_{m n}
 \kappa^{m} \lambda^{n}
\end{equation}
by identifying $\mathbb{K}^\uparrow=\text{Hom}_\mathbb{K}(H,\mathbb{K})$ with
$\cal H$ via the pairing.
The equivariance condition (\ref{equivariancecondition})
\begin{equation}
   \langle hl, f\rangle = \langle h, f \rangle \dashv l, \qquad
    \forall l \in L, \quad \forall h \in H,
\end{equation}
together with duality  give  the relation
\begin{equation}
     m! n! f_{m n}= \langle k_m l_n, f \rangle =
     \langle k_m, f \rangle  a_n = m! f_{m 0} a_n.
\end{equation}
Hence, the elements of $\mathbb{K}^\uparrow$ are
\begin{equation}
f= \kappa \psi, \qquad \kappa \in {\cal K},\ \ 
\psi=e^{a_1 \lambda_1} \cdots e^{a_s \lambda_s} ,
\end{equation}
where  $\cal K$ is the  subspace of $\cal H$ generated by the linear  combinations
of the ordered monomials $(\kappa^m)_{m \in \mathbb{N}^r}$. Because 
$\psi$ is product of exponentials,
 ${\cal K}$ and $\mathbb{K}^\uparrow$ are  isomorphic 
($\kappa  \to  \kappa \psi$).

The action of $H$ on  $\mathbb{K}^\uparrow$ is determined
knowing the action over the basis elements $(\kappa^p \psi)_{p \in \mathbb{N}^r}$. Putting
\begin{equation}
(\kappa^p \psi) \dashv h =
   \sum_{(m,n)\in \mathbb{N}^r \times \mathbb{N}^s}  [h]^p_{m n} \kappa^m\lambda^n, 
\qquad  p \in \mathbb{N}^r,
\end{equation}
the coefficients $[h]^p_{m n}$ are evaluated by means of duality
 \begin{equation}\label{coef} 
m! n! [h]^p_{m n} = \langle (\kappa^p \psi) \dashv h,
 k_m l_n\rangle =
    \langle \kappa^p \psi,  h k_m l_n\rangle =
    \langle \kappa^p \psi,  h k_m \rangle a_n.
 \end{equation}

The properties of the action allow to compute it only for the generators of $H$ instead 
of considering an arbitrary element of $H$. The problem  reduces to write $h k_m$ in
normal ordering to get the value of the paring in (\ref{coef}). Since in many cases this
task is very cumbersome, our objective will be to take advantage of the bicrossproduct
structure to simplify the computations.

In the following we will restrict ourselves to Hopf algebras having a
bicrossproduct structure like  $H= {\cal K} \bicross {\cal L}$, such 
that ${\cal K}$ is cocommutative and ${\cal L}$ commutative. Let us suppose that the 
algebras  $\cal K$ and $\cal L$  are finite  generated by
the sets $\{ k_i \}_{i=1}^r$ and  
$\{ l_i \}_{i=1}^s$, respectively,  the $k_i$'s are primitive and  $(k_n)_{n\in
\mathbb{N}^r}$ and $(l_m)_{m\in \mathbb{N}^s}$ are bases of the vector spaces underlying   to
$\cal K$ and $\cal L$, respectively.
Let ${\cal K}^*$ and ${\cal L}^*$ be the dual algebras of $\cal K$ and $\cal L$  having  dual
systems to those of $\cal K$ and $\cal L$ with analogue properties  to them. So,  duality
between  $H$ and $H^*$ is given  by  (\ref{pairing}).

We are interested in the construction of the representations induced by 
`real' characters of the commutative sector  $\cal L$.
We will show that the solution of this problem can be reduced to the study of
certain dynamical systems which  present, in general, a nonlinear action.
The above results can be summarized in the following theorem.
\begin{theorem}\label{tind0} 
The   carrier space
of  the  representation of $H$ induced by the  character $a$ of  $\cal L$ 
\begin{equation}  \label{caracter}
 1 \dashv l_n= a_n, \qquad a \in \mathbb{C}^s , \ \ n \in \mathbb{N}^s,\ 
\end{equation}
is isomorphic to ${\cal K}^*$ and is constituted by the  elements of  the  form  
\begin{equation}
 \kappa \psi, \ \ \ \kappa \in {\cal K}^* ,\qquad       
\psi= e^{a_1 \lambda_1} e^{a_2 \lambda_2} \cdots  e^{a_s \lambda_s}, \ \ \
\lambda_i \in {\cal L}^* .
\end{equation}
The  induced action  of the elements of $H$ over the elements of the carrier space
$\mathbb{C}^\uparrow$ is given by 
\begin{equation} \label{accion0}
           f \dashv h= \sum_{m\in \mathbb{N}^r}
            \kappa^m \langle h \frac{k_m}{m!}, f \rangle \psi,
            \qquad  f\in \mathbb{C}^\uparrow.
\end{equation}
\end{theorem}

To obtain the  explicit action  of the generators of $\cal K$ and $\cal L$ in  the
induced  representation, let us start  identifying  $\cal L$, because 
is commutative, with the algebra of functions $F(\mathbb{R}^s)$ by the  morphism
${\cal L} \to F(\mathbb{R}^s),\ (l \mapsto \tilde{l})$,
 mapping the  generators of  $\cal L$ into the canonical projections
\begin{equation}
 \tilde{l}_j(x)= x_j, \qquad\quad 
\forall x=  (x_1, x_2, \ldots, x_s)\in \mathbb{R}^s, \quad
 j=1,2,\ldots, s.
\end{equation}
This identification allows to chose  a  $*$--structure keeping
invariant the  generators of $\cal L$ by
\begin{equation} \label{estre}
       l_j^*= l_j, \qquad j=1,2,\ldots, s.
\end{equation}
Hence, the characters of $\cal L$ (\ref{caracter}) compatible with (\ref{estre}) are real. 
They can be written now as 
\begin{equation}
  1 \dashv l = \tilde{l}(a), \qquad a \in \mathbb{R}^s.
\end{equation}
Also, the  right action  of $\cal K$ on $\cal L$ can be carried to 
 $F(\mathbb{R}^s)$.  Since the generators of $\cal K$ are primitive, they
act by derivations on the   $\cal K$--module algebra of ${\cal K} \RL {\cal L}$ inducing vector
fields, $X_i$,  on $\mathbb{R}^s$ by 
\begin{equation}
 X_i\, \tilde{l}= \widetilde{l\ract k_i},
   \qquad i=1,2, \ldots, r.
\end{equation}
The  associated flow to $X_i$,
 $\Phi_i:\mathbb{R} \times \mathbb{R}^s \rightarrow \mathbb{R}^s$,
is given by
\begin{equation} \label{flujocampo}
   (X_i f)(x)= (D f_{x,\Phi_i})(0),
\end{equation}
where $f_{x,\Phi_i}(t)= f\circ \Phi_i^t(x)$ and $D$ is the  derivative  operator. Thus, we can
state the following theorem that allows us to have explicit formulae for the action.
\begin{theorem}\label{tind1}
The  explicit action  of the generators of  $\cal K$ and $\cal L$ in
 the  induced representation   determined in Theorem~\ref{tind0} and
realized in the   space ${\cal K}^*$ is given by the following  expressions:
\begin{equation} \label{accion1}
\begin{array}{rcl}
           \kappa \dashv k_i & = & \kappa \prec k_i \  ,\\[0.2cm]
           \kappa \dashv l_j & = & \kappa \, \,
        \hat{l}_j  \! \circ \Phi_{(\kappa_1,\kappa_2, \ldots, \kappa_r)}(a) \ ,
\end{array}
\end{equation}
where $i \in \{1,\ldots, r\}$, $j \in \{1,\ldots, s\}$, the   symbol
$\prec$ denotes the  regular action of  $\cal K$ on ${\cal K}^*$ and
$\Phi_{(\kappa_1,\kappa_2, \ldots, \kappa_r)}=
\Phi_r^{\kappa_r} \circ \cdots \circ   \Phi_2^{\kappa_2} \circ 
\Phi_1^{\kappa_1}$.
\end{theorem}

We can reformulate the induction procedure in terms of modules because their
deep relation with  representations \cite{olmo01,olmo00}.
The regular $H$--modules associated to a
Hopf algebra  $H$ and its dual $H^*$,  
$(H,\prec, H)$, $(H^*,\succ, H)$, $(H,\succ, H)$ and
$(H^*,\prec,H)$, will help us to
describe the actions as well as the carrier spaces involved in the induced representations.
The modules  $(H^*,\succ, H)$ and
$(H^*,\prec,H)$ are called co-spaces since they can be considered as an algebraic
generalization of the concept of $G$--space.

The commutative or cocommutative Hopf algebras are, in fact,  of the form
$F(G)$ or $\mathbb{C}[G]$, the group algebra of $G$,  (or $U(\mathfrak{g})$) for any  group
$G$ \cite{Maj95a}. Hence, the  bicrossproduct Hopf algebras, that we consider in
this paper, can  be factorized as  
 \begin{equation}
     H= U(\mathfrak{k})  \, \RL F(L),\qquad  H^*= F(K) \LR  \, U(\mathfrak{l}).
 \end{equation}
where  $K$ and $L$
are Lie groups with associated  Lie algebras $\mathfrak{k}$ and
$\mathfrak{l}$, respectively. 

The use of  elements of  $H$ and $H^*$, like 
\begin{equation}\begin{array}{lllll}  \label{tablita}
&  k \lambda \in H , &  \qquad\qquad  & k \in K ,\ \  &\lambda \in F(L) ,\\[0.2cm]
& \kappa l \in H^* ,  &\qquad\qquad &\kappa \in F(K) , \ \  & l \in L ,
\end{array}\end{equation}
instead of the standard bases of  ordered  monomials, allows an effective  description  of  the
$H$--regular modules. 
\begin{theorem}\label{crbicross}
The  action  on each of  the four regular $H$--modules is:
\begin{equation}\begin{array}{lllll} 
(H,\prec, H)&:&\ \ (k\lambda) \prec k'=  kk'(\lambda \ract k') , 
\qquad\quad  & (k\lambda) \prec \lambda'=  k\lambda \lambda' ;\\[0.2cm]
 (H^*,\succ, H)&:&\ \  k' \succ (\kappa l)  =  (k'\succ \kappa) (k'\lact l),
\qquad\quad  &\lambda' \succ (\kappa l)  = \lambda'( l)\kappa  l ;\\[0.2cm]
(H,\succ, H)&:&\ \  k' \succ (k\lambda) =  k'k\lambda ,
\qquad\quad  &\lambda' \succ (k\lambda) =  k (\lambda' \ract k) \lambda ;\\[0.2cm]
(H^*,\prec, H)&:&\ \   (\kappa l) \prec k'  =  (\kappa \prec k') l ,
\qquad\quad  &(\kappa l) \prec \lambda'  =\kappa (\lambda'\circ \hat{l}) l ;
\end{array}
\end{equation}
where  $\hat{l}$ is the  map, 
$ K  \rightarrow L \ (k  \mapsto  k \lact l)$, projecting the  group $K$ on the  orbit passing
through $l \in L$ such that  $\langle l^{(1)}, \lambda \rangle l^{(2)} = \lambda \circ
\hat{l}$, for all $k,k' \in K$, $\lambda, \lambda' \in F(L)$, $\kappa \in F(K)$ and $l \in L$.
\end{theorem} 
Notice that in this  theorem we have taken into account that in
$(U(\mathfrak{l}),\succ, F(L))$ holds 
\begin{equation}
 \lambda \succ l= \lambda(l) l,
  \qquad \forall \lambda \in F(L), \; \forall l \in L .
\end{equation}

Now we can carry a
complete  analysis of the  representations of $H= U(\mathfrak{k}) \RL F(L)$ induced by the
one-dimensional modules of  its  commutative sector. Note that the set of  characters
of the algebra $F(L)$ is its spectrum and, hence,  the spectrum of $F(L)$ is isomorphic to $L$.
 Fixed an element $l$ of $ L$, the  character is given by
\begin{equation}  \label{car}
   1 \dashv \lambda = \lambda(l), \quad \lambda \in F(L).
\end{equation}
The carrier space 
$\mathbb{C}^\uparrow$ of the  representation of  $H= U(\mathfrak{k}) \RL F(L)$
induced by (\ref{car}) is the set of   elements $f \in H^*$ satisfying  the 
condition
 \begin{equation}
    \lambda \succ f = \lambda(l) f, \qquad \forall \lambda \in F(L).
 \end{equation}
Expanding $f$ in terms of the  bases  of $\mathfrak{k}$ and $\mathfrak{l}$, imposing the
equivariance condition and considering the  definition of 
second kind coordinates $\lambda_j$ over the   group $L$ we obtain that 
\begin{equation}
     f= \kappa l, \qquad  \kappa \in F(K).
\end{equation}

The  right regular action  determines  
the  action  on the  induced module that can be carried to
$F(K)$ using the isomorphism
$ F(K)  \rightarrow  \mathbb{C}^\uparrow \ (  \kappa  \mapsto  \kappa l)$. So,
\begin{equation}
 \kappa \dashv k= \kappa \prec k, \qquad
 \kappa  \dashv \lambda=   \kappa  (\lambda \circ \hat{l}).
\end{equation}
Comparing these expressions with those of   Theorem~\ref{tind1} we see 
that  the  action of  $U(\mathfrak{k})$ is determined by 
the  regular action.  The  action of   $F(L)$ is 
 multiplicative and  its evaluation  is essentially reduced to  obtain  the 
flows associated to  the  action of  $K$ on $L$ derived from  the bicrossproduct 
structure  of  $H$.

Finally, next theorem summarizes the induction procedure for bicrossproduct algebras.
\begin{theorem}  \label{repequiv}
Let us consider an element  $l\in L$ supporting a global  action of the  group
$K$.  The carrier space, $\mathbb{C}^\uparrow$, of the  representation of  $H$
induced by the  character    determined by $l$ is the  set of  elements of  $H^*$
of  the  form
 \begin{equation} \kappa l, \qquad \kappa \in F(K).
\end{equation}
There is an isomorphism between
$\mathbb{C}^\uparrow$ and $F(K)$ given by the  map $\kappa \mapsto \kappa l$.
The  action  induced  by the elements of  the  form 
$k \in K$ and $\lambda \in F(L)$ in the  space $F(K)$ is
\begin{equation}
\begin{array}{lll}
        \kappa \dashv k &= & \kappa \prec k \\[0.2cm]
        \kappa  \dashv \lambda &= &
             \kappa  (\lambda \circ \hat{l}) .
\end{array}
\end{equation} 
The modules induced by $l$ and $k\lact l$ are isomorphic. So,
the induction algorithm establishes a correspondence between the  
space of  orbits $L/K$ and the set of  equivalence classes of  representations.
\end{theorem}
\subsect{Local  representations}
\medskip

The quantum counterpart of the  local representations \cite{olmo84} of Lie groups can be
obtained  inducing from representations of  the  subalgebra
$U(\mathfrak{k})$. 
Given a character $\kappa$ of  $U(\mathfrak{k})$, 
  \begin{equation} 
\kappa \in \text{Spectrum} \, U(\mathfrak{k}) \subset F(K), \qquad
      k \vdash 1 = \kappa(k),
\end{equation}
the  carrier space, $\mathbb{C}^\uparrow$, of  the  representation induced by $\kappa$ is
determined by the equivariance condition
 \begin{equation}
    f \prec k=  \kappa(k) f, \qquad  \forall k \in U(\mathfrak{k})  
 \end{equation} 
obtaining that the elements of  $\mathbb{C}^\uparrow$ are of the form
\begin{equation}
         f = \kappa l, \qquad  l \in U(\mathfrak{l}).
\end{equation}
 The  isomorphism $U(\mathfrak{l})\to \mathbb{C}^\uparrow$ ($l \mapsto \kappa l$),
allows to realize the  induced representation    over $U(\mathfrak{l})$ 
\begin{equation} \label{reploc}
       k \vdash l= \kappa(k) \; k \lact l ,\qquad
        \lambda \vdash l = \lambda(l) l .
\end{equation}

\sect{Quantum  $\mathfrak{iso}_{\omega}(2)$ algebras}
\label{cayley-kleinalgebras}
 
The called  CK  pseudo-orthogonal algebras $\mathfrak{so}_{\omega_1,\omega_2, \ldots,
\omega_N}(N+1)$ is a family of ${(N+1)N}/{2}$ dimensional real Lie  algebras 
characterized by $N$ real parameters $(\omega_1,\omega_2, \ldots, \omega_N)$
\cite{olmo93,olmo97}. In the `geometric' basis $(J_{ij})_{0\leq i < j \leq N}$ the
nonvanishing commutators are
\begin{equation} 
[J_{ij}, J_{ik}]= \omega_{ij}J_{jk}\qquad
   [J_{ij}, J_{jk}]= - J_{ik} \qquad
   [J_{ik}, J_{jk}]= \omega_{jk}J_{ij},
\end{equation}
with  
$0<i<j<k<N$ and $\omega_{ij}={\displaystyle \prod_{s=i+1}^j \omega_s}$.
 
The parameters $\omega_i$, in fact, only take  the values $1,\ 0$ and $-1$ since the generators
$J_{ij}$ can be rescaled.  When  $\omega_i\neq 0 \ \forall i$, the Lie
algebra $\mathfrak{so}_{\omega_1,\omega_2, \ldots, \omega_N}(N+1)$
is isomorphic to  some of the  pseudo-orthogonal algebras
$\mathfrak{so}(p,q)$ with $p+q=N+1$ and $p\geq q \geq 0$. If some of the 
$\omega_i$ vanish  the corresponding  algebra is inhomogeneous. 

In this paper we are interested in the particular case of $N=2$ and $\omega_1=0$. These
inhomogeneous  algebras
$\mathfrak{so}_{0,\omega_2}(3)$ can be realized as algebras of  groups
of  affine transformations on $\mathbb{R}^2$ \cite{olmo93} (i.e, the Euclidean group when
$\omega_2=1$, the Galilei group  when $\omega_2=0$ and the Poincar\'e group when
$\omega_2=-1$).  In this case the generators
$J_{0i}$ are denoted by $P_i$ stressing  their role as 
generators of  translations. The   remaining generator $J_{12}\equiv J$ is associated to 
compact ($\omega_2=1$) and noncompact rotations (galilean  and
lorentzian boosts).

The algebra $\mathfrak{so}_{0,\omega_2}(3)$, that we shall  denote by
$\mathfrak{iso}_{\omega}(2)$, is characterized by
the  commutators
\begin{equation} 
    [J, P_1]= P_2, \qquad [J, P_2]= - \omega P_1, \qquad  [P_1, P_2]= 0 .  
\end{equation}

We shall use the following `generalized trigonometric  
 functions' \cite{olmo93,olmo93a}
\begin{equation} \label{scgen}
  {\rm C}_{\omega}(x)= \frac{e^{\sqrt{-\omega} x} + e^{-\sqrt{-\omega} x}}{2},
 \qquad
  {\rm S}_{\omega}(x)= \frac{e^{\sqrt{-\omega} x} - e^{-\sqrt{-\omega} x}}{2\sqrt{-\omega}}.
\end{equation}
When $\omega <0$ ($\omega >0$) these expressions become
the trigonometric  (hyperbolic) functions. 
For  $\omega=0$ the parabolic functions ${\rm C}_0(x)=1$ and ${\rm S}_0(x)=x$ are obtained.
The functions (\ref{scgen}) satisfy identities similar to those of the usual 
trigonometric functions. Some useful properties that  will be used in next computations
are:
\begin{equation}
  \begin{array}{c}
  {\rm C}_\omega^2(x) + \omega {\rm S}_\omega^2(x)=1, \\[0.2cm]
   {\rm C}_\omega(x+y)=   {\rm C}_\omega(x)   {\rm C}_\omega(y) -  
\omega  {\rm S}_\omega(x)  {\rm S}_\omega(y), \quad
  {\rm S}_\omega(x+y)=   {\rm S}_\omega(x)   {\rm C}_\omega(y) +  
{\rm C}_\omega(x)   {\rm S}_\omega(y), \\[0.2cm]
  C'_{\omega}(x)= - \omega {\rm S}_\omega(x), \qquad
    S'_{\omega}(x)=  {\rm C}_\omega(x).
  \end{array}
\end{equation} 

A  simultaneous standard deformation for all the  enveloping
algebras  $U(\mathfrak{so}_{\omega_1, \omega_2}(3))$ was introduced in \cite{olmo93a},
 The particular case of the deformed Hopf algebras $U_\lambda(\mathfrak{iso}_{\omega}(2))$
is obtained taking $\omega_1=0$ and, obviously, $\omega_2=\omega$.

 It was proved  in \cite{olmo97}  that the standard quantum Hopf algebras
$U_\lambda(\mathfrak{iso}_{\omega_2, \omega_3,\ldots, \omega_N }(N))$
have a structure of bicrossproduct.  So,
$U_\lambda(\mathfrak{iso}_{\omega}(2))$ is characterized in  a basis adapted to its
bicrossproduct structure  by 
 $$
\begin{array}{c}
  [J,P_1]= \frac{1- e^{-2 \lambda P_2}}{2 \lambda} + \frac{1}{2} \lambda \omega P_1^2, 
\qquad [J,P_2]= - \omega P_1;   \\[0.4cm]
  \Delta(P_1)= P_1 \otimes 1 + e^{- \lambda P_2} \otimes  P_1,  \quad
   \Delta(P_2)= P_2 \otimes 1 + 1 \otimes P_2, \quad
   \Delta(J)= J \otimes 1 + e^{- \lambda P_2} \otimes J;  \\[0.4cm]
  \epsilon(P_1)= \epsilon(P_2)= \epsilon(J)=0; \\[0.4cm] 
    S(P_1)= -e^{\lambda P_2} P_1, \quad S(P_2)=-P_2, \quad
    S(J)= -e^{\lambda P_2} J.
 \end{array}
$$

The  bicrossproduct structure 
$U_\lambda (\mathfrak{iso}_{\omega}(2))=  {\cal K} \RL {\cal L} $, where  $\cal L$ is  the 
Hopf subalgebra spanned by  ($P_1$, $P_2$) and $\cal K$ is the commutative and cocommutative
Hopf algebra  generated by  $J$, is determined by  the  action of 
$\cal K$ on $\cal L$ given by
\begin{equation}
  P_1 \ract J=[P_1, J]=
     - \left[ \frac{1- e^{-2 \lambda P_2}}{2 \lambda} +
     \frac{1}{2} \lambda \omega P_1^2\right], \qquad
  P_2 \ract J=[P_2, J]= \omega P_1,
\end{equation} 
and  the  left coaction of  $\cal L$ on $\cal K$, which over the 
generator of  $\cal K$ takes the value
\begin{equation}
   J \lcact = e^{-\lambda P_2} \otimes J.
\end{equation}

The dual algebra $F_\lambda(ISO_{\omega}(2))$
is generated by the  local coordinates
$\varphi, a_1, a_2$.  Its
commutators, coproduct, counit  and antipode are given by
\begin{equation}
 \begin{array}{c}
  [a_1, \varphi]= \lambda(1- {\rm C}_\omega(\varphi)), \qquad
   [a_2, \varphi]= \lambda {\rm S}_\omega(\varphi), \qquad
   [a_1, a_2]= \lambda a_1;     \\[0.4cm]
  \Delta(a_1)= a_1 \otimes {\rm C}_\omega(\varphi) + a_2 \otimes \omega {\rm S}_\omega(\varphi)
    + 1 \otimes a_1, \\[0.4cm]
   \Delta(a_2)=- a_1 \otimes {\rm S}_\omega(\varphi) + a_2 \otimes  {\rm C}_\omega(\varphi)
 + 1 \otimes a_2, \\[0.4cm]
   \Delta(\varphi)= \varphi \otimes 1 + 1 \otimes \varphi;  \\[0.4cm]
   \epsilon(a_1)=  \epsilon(a_2)= \epsilon(\varphi)= 0; \\[0.4cm]
   S(a_1)= - {\rm C}_\omega(\varphi) a_1 - \omega {\rm S}_\omega(\varphi) a_2, \quad
   S(a_2)=  {\rm S}_\omega(\varphi) a_1 - {\rm C}_\omega(\varphi) a_2, \quad
   S(\varphi)= -\varphi.
\end{array}
\end{equation}
This Hopf algebra exhibits  the  bicrossproduct structure
 $F_\lambda(ISO_{\omega}(2))= {\cal K}^* \LR {\cal L}^*$,
 dual of  the  above one, with
${\cal K}^*$ generated by $\varphi$ and ${\cal L}^*$ by $a_1, a_2$.
The  left action of  ${\cal L}^*$ over ${\cal K}^*$ is given by 
\begin{equation}
   a_1 \lact \varphi= \lambda(1- {\rm C}_\omega(\varphi)), \qquad
   a_2 \lact \varphi= \lambda {\rm S}_\omega(\varphi).
\end{equation}
The  right coaction  of  ${\cal K}^*$ over  ${\cal L}^*$ takes the following values
over  the generators of  ${\cal L}^*$
\begin{equation}
 \rcact a_1= a_1 \otimes {\rm C}_\omega(\varphi) +  a_2 \otimes \omega {\rm S}_\omega(\varphi), \qquad
 \rcact a_2= - a_1 \otimes  {\rm S}_\omega(\varphi) +  a_2 \otimes  {\rm C}_\omega(\varphi).
\end{equation} 

The results mentioned in the last paragraph of Section~\ref{preliminaries} about 
dual bases and
$*$--structures over bicrossproduct Hopf algebras allow  to construct a pair of  
dual bases in such a way that
the duality form between  
$U_\lambda (\mathfrak{iso}_{\omega}(2))$ and
$F_\lambda(ISO_{\omega}(2))$ is
\begin{equation}
  \langle J^m P_1^n P_2^p, \varphi^q a_1^r a_2^s\rangle= m! n! p! \;
       \delta_m^q \delta_n^r \delta_p^s.
\end{equation}

\sect{Representations of  $U_\lambda(\mathfrak{iso}_{\omega}(2))$}
\label{representationsquantumstandardulambda}

Firstly, according  to the theory of representations of bicrossproduct algebras displayed in
Section~\ref{inducedrepresentations} and, in particular, Theorem~\ref{tind1}, we need to know
the flow associated to the action of $\cal K$ on $\cal L$.


\subsect{One-parameter flow for $U_\lambda(\mathfrak{iso}_{\omega}(2))$}
\label{oneparameterflow}

In Ref.~\cite{olmo97} the factor  $\cal K$ is interpreted as the  enveloping algebra 
$U(\mathfrak{so_\omega(2)})$, while $\cal L$ is seen as a deformation
of the algebra of the  group of  the translations in the plane $T_2$,
$U_{\lambda}(\mathfrak{t}_2)$. Hence,
  \begin{equation}
     U_{\lambda}(\mathfrak{iso_\omega(2)}) =
          U(\mathfrak{so_\omega(2)}) \RL  U_{\lambda}(\mathfrak{t}_2), \qquad
F_\lambda(ISO_{\omega}(2))=F(SO_{\omega}(2)) \LR F_\lambda(T_2). 
\end{equation} 
However, the interpretation given here by us is different since 
  \begin{equation}
       U_\lambda(\mathfrak{iso_\omega(2)}) =
       U(\mathfrak{so_\omega(2)}) \RL  F(T_{\lambda,2}),
  \end{equation}        
being $T_{\lambda,2}$ the  Lie group with  composition law
 \begin{equation}
    (\alpha'_1, \alpha'_2)(\alpha_1, \alpha_2)=
    (\alpha'_1+ e^{- \lambda \alpha'_2} \alpha_1, \alpha'_2+ \alpha_2).
 \end{equation} 
Now the generators $P_1$ and $P_2$ are considered as a global coordinate system over 
$T_{\lambda,2}$, i.e., 
 \begin{equation}
 P_1(\alpha_1, \alpha_2)= \alpha_1, \qquad
           P_2(\alpha_1, \alpha_2)= \alpha_2.
\end{equation}
The  module algebra structure included in  
$U_\lambda(\mathfrak{iso_\omega(2)}) =  U(\mathfrak{so_\omega(2)}) \RL  F(T_{\lambda,2})$
determines  the  action of  $SO_\omega(2)$ over  $T_{\lambda,2}$ by means of the  vector
field 
\begin{equation}
    \hat{J}_{\omega, \lambda} = -\left[\frac{1- e^{-2 \lambda P_2}}{2 \lambda} +
  \frac{1}{2} \lambda \omega P_1^2\right] \frac{\partial}{\partial P_1} +
              \omega P_1 \frac{\partial}{\partial P_2}.
\end{equation}

Prior to determine the flow associated to $\hat{J}_{\omega, \lambda}$ we can obtain 
a first qualitative  information  by identifying the fixed points of
$\hat J_{\omega,\lambda}$. When $\omega\not = 0$ only the origin $(0,0)$ is an   equilibrium
point, but if $\omega=0$  the   set of  fixed points is the straight line  of  equation
$\alpha_2=0$.  The  deformation due to $\lambda$ does not change the  character of these
points. Thus, the   point is elliptic if $\omega> 0$ and   hyperbolic if $\omega <0$.

The  one-form
\begin{equation}
 \eta_\rho= \rho\left[ \omega P_1 dP_1 + (\frac{1-e^{-2\lambda P_2}}{2\lambda} +
 \frac{1}{2}\lambda \omega P_1^2) dP_2\right], \qquad  \rho \in F(T_{\lambda,2}),
\end{equation}
verifies $\hat{J}_{\omega,\lambda} \rfloor \eta_\rho=0$. Taking $\rho_0= \lambda^2 e^{\lambda
P_2}$ one obtains that $\eta_{\rho_0}$  is  exact and invariant under
$\hat{J}_{\omega,\lambda}$. Adding  a constant and rescaling the invariant function under  the 
action of  $J$
    \begin{equation}\label{casimirh}
     h= \frac{1}{2} \omega \lambda^2 P_1^2 e^{\lambda P_2}+ \cosh(\lambda P_2),
    \end{equation}
we obtain a central element of $U_\lambda (\mathfrak{iso}_{\omega}(2))$
\begin{equation}\label{casimir}
 {\rm C}_{\omega,\lambda}= \omega  P_1^2 e^{\lambda P_2}+
2 \frac{ \cosh(\lambda P_2)-1 }{\lambda^2}=  \omega  P_1^2 e^{\lambda P_2}+
      \frac{4}{\lambda ^2} \sinh^2(\frac{\lambda}{2} P_2) 
 \end{equation} 
 such  that in the  limit $\lambda \to 0$ we recover the nondeformed Casimir 
${\rm C}_{\omega, 0}= \omega P_1^2 + P_2^2$ of $\mathfrak{iso}_{\omega}(2)$.

The   computation of the trajectories of
$\hat{J}_{\omega,\lambda}$ requires to solve the equation system
\begin{equation} \label{ecumoviso}
 \begin{split}
     \dot{\alpha_1}= & - \left[ \frac{1- e^{-2\lambda \alpha_2}}{2 \lambda}+
       \frac{1}{2} \lambda \omega \alpha_1^2\right],  \\[0.3cm]
       \dot{\alpha_2}= &\ \omega \alpha_1.
    \end{split}
 \end{equation}
The  integral curve $\gamma$ can be expressed as
 \begin{equation}\begin{split}
        \alpha_1(t)= &  - \frac{1}{\lambda} \frac{\sinh(\lambda \beta) 
{\rm S}_\omega(t)}{  \cosh(\lambda \beta) + \sinh(\lambda \beta) {\rm C}_\omega(t)},\\[0.3cm]
 \alpha_2(t)= & \frac{1}{\lambda} \ln[\cosh(\lambda \beta) 
+ \sinh(\lambda \beta) {\rm C}_\omega(t)].
\end{split}\end{equation} 
Using the fact that $\Phi^t_{\omega, \lambda}(\gamma(\tau))= \gamma(\tau+t)$ the   flow can be  
evaluated
\begin{equation} \label{flowiso}
\begin{split}
\Phi_{\omega,\lambda} ^t &( \alpha_1,\alpha_2)=\\[0.3cm]
   & \left( \frac{2\alpha_1 \lambda e^{\lambda \alpha_2} {\rm C}_\omega (t)
 + (\omega \lambda ^2 \alpha_1^2 e^{\lambda \alpha_2} 
- {\sinh(\lambda\alpha_2)}) {\rm S}_\omega(t)}{2\lambda \cosh(\lambda \alpha_2) 
+  \omega \lambda^3 \alpha_1^2 e^{\lambda \alpha_2}+
    (2\lambda \sinh(\lambda \alpha_2) -   \omega \lambda^3 \alpha_1^2 
e^{\lambda \alpha_2}) {\rm C}_\omega(t) +
  2\omega \lambda ^2 \alpha_1 e^{\lambda \alpha_2} {\rm S}_\omega(t)}\right. , \\[0.25cm]
  &   \frac{1}{ \lambda} \ln\left[\cosh(\lambda \alpha_2) +
 \left.   \frac{1}{2} \omega \lambda^2 \alpha_1^2 e^{\lambda \alpha_2}
+  (\sinh(\lambda \alpha_2) -
  \frac{1}{2} \omega \lambda^2 \alpha_1^2 e^{\lambda \alpha_2})  {\rm C}_\omega (t) +
\omega \lambda \alpha_1 e^{\lambda \alpha_2} {\rm S}_\omega (t)\right]\right).
 \end{split}
\end{equation}
In the   limit $\lambda \to 0$  we recover the   linear flow
  \begin{equation}
 \Phi_{\omega,0}^t(\alpha_1,\alpha_2)=
     \left( {\rm C}_\omega(t) \alpha_1 - {\rm S}_\omega(t) \alpha_2, \ 
         \omega {\rm S}_\omega(t) \alpha_1 + {\rm C}_\omega(t) \alpha_2\right),
  \end{equation}       
 which corresponds to the nondeformed action given by the vector field 
$\hat{J}_{\omega,0}= - P_2 {\partial}_{P_1} + \omega P_1 {\partial}_{P_2}$. 

Note that  although  expression (\ref{flowiso})  has been  obtained for $\omega >0$ and
$\lambda >0$ it can be proved that it is also valid for the   remaining values. In
some sense, the flow has an ``analytic'' dependence on the parameters $\omega$ and
$\lambda$ which allows to extend   the results obtained in a region of the parameters space
to whole it.

The   flow (\ref{flowiso}) is globally defined when $\omega \geq 0$ but is only local for
$\omega <0$ as it is easy to prove considering the integral curve passing through a  point
like $(0,\alpha_2)$: since the logarithm argument has to be positive 
one gets the  inequality
 \begin{equation} 
\cosh(\lambda \alpha_2) + \sinh(\lambda \alpha_2) {\rm C}_\omega(t) > 0 ,
\end{equation} 
if the   product $\lambda \alpha$ is negative the definition interval of 
$t$ is bounded
\begin{equation}
 t \in (- {\rm C}_\omega^{-1}(-\coth(\lambda \alpha_2)),\
      {\rm C}_\omega^{-1}(-\coth(\lambda \alpha_2)) ). 
\end{equation}

The   flow (\ref{flowiso}) describes  the  action of  $SO_\omega(2)$ over  $T_{\lambda,2}$
 \begin{equation} 
e^{tJ}\lact (\alpha_1, \alpha_2)=  \Phi^t_{\lambda, \omega}(\alpha_1, \alpha_2),
\end{equation}
which  decomposes the space $T_{\lambda,2}$
on strata of orbits  depending  on the values of $\omega$:

\noindent $\bullet$ Case $\omega<0$

- A   stratum with only one orbit with only one point: (0,0).
The isotropy group is $SO_\omega(2)$.
       
- Four orbits determined by the points
\begin{equation}
   ( \frac{1- e^{- \lambda }}{\lambda \sqrt{-\omega}},  1), \quad
   ( \frac{1- e^{ \lambda }}{\lambda \sqrt{-\omega}},  -1), \quad
   (- \frac{1- e^{- \lambda }}{\lambda \sqrt{-\omega}},  1), \quad
   (- \frac{1- e^{ \lambda }}{\lambda \sqrt{-\omega}},  -1).
\end{equation}

These orbits are isomorphic to  $\mathbb{R}$ and have the point $(0,0)$ as accumulation point.

- The   remaining points of  $T_{\lambda,2}$ constitute another stratum fibrered by
orbits diffeomorphic  to  
$\mathbb{R}$ since they are branches of deformed hyperbolae.

\noindent $\bullet$  Case $\omega=0$
 
- A stratum is constituted by  the points $(\alpha_1,0)$, each of  them is an orbit.
        
- The  orbits 
$\{ (\alpha_1, \alpha_2) \; | \; \alpha_1 \in \mathbb{R},\ \alpha_2\neq 0\  {\rm but \ fixed}\}$
determine a stratum.

\noindent  $\bullet$ Case $\omega>0$

- The   point   $(0,0)$ constitutes the only orbit of this stratum.
        
- The  orbits diffeomorphic to the  circle determine another stratum.
\medskip

Note that the deformation associated to $\lambda$ does not give qualitative changes with respect
to  the  nondeformed case. Summarizing, we can say that the quotient spaces are
isomorphic
 \begin{equation}
       T_{\lambda,2}/ SO_\omega(2) \simeq  T_{0,2}/ SO_\omega(2).
 \end{equation}


\subsect{Regular co-spaces}\label{regularcospaces}

Once the flow associated to the bicrossproduct structure of
$U_\lambda(\mathfrak{iso}_{\omega}(2)))$ is known, we can study the regular  co-spaces 
$(F_\lambda(ISO_{\omega}(2)), \prec ,U_\lambda(\mathfrak{iso}_{\omega}(2)))$ and 
$(F_\lambda(ISO_{\omega}(2)), \succ ,U_\lambda(\mathfrak{iso}_{\omega}(2)))$,  which are basic
elements to characterize the induced representations of
$U_\lambda(\mathfrak{iso}_{\omega}(2)))$.

The  structures of both  co-spaces are easily deduced combining  Theorem~\ref{crbicross} with 
expression (\ref{flowiso}) of the flow of $\hat{J}$. 
Firstly, remember that $(F_\lambda(ISO_{\omega}(2))$ can be described by elements of
the form
\begin{equation}
    \phi(\alpha_1,\alpha_2),\qquad\qquad \phi \in F(SO_{\omega}(2)), \quad
       (\alpha_1,\alpha_2) \in T_{{\lambda},2},
\end{equation}
instead of  monomials $\varphi^q a_1^r a_2^s$.
\bigskip

{\bf 1).-} For the right coregular module $(F_\lambda(ISO_{\omega}(2)), \prec
,U_\lambda(\mathfrak{iso}_{\omega}(2)))$ we have that 
$$
\begin{array}{l}
    [\phi(\alpha_1,\alpha_2)]\prec e^{tJ}=
         \phi(e^{tJ} \, \cdot \,)(\alpha_1,\alpha_2), \\[0.3cm]
    [\phi(\alpha_1,\alpha_2)]\prec P_1= \\[0.15cm]
 \phi {\displaystyle\frac{ 2 \alpha_1 \lambda e^{\lambda \alpha_2}
{\rm C}_\omega(\varphi) +
     ( \omega \lambda^2 \alpha_1^2 e^{\lambda \alpha_2}
        -{2\sinh(\lambda \alpha_2)})  {\rm S}_\omega(\varphi)}
{ 2 \lambda \cosh(\lambda \alpha_2) +  \omega \lambda^3 \alpha_1^2
 e^{\lambda  \alpha_2})+ (2\lambda \sinh(\lambda \alpha_2) -
     \omega \lambda^3 \alpha_1^2 e^{\lambda \alpha_2})  {\rm C}_\omega(\varphi) +
       2\omega \lambda^2 \alpha_1 e^{\lambda \alpha_2} {\rm S}_\omega(\varphi)}}
(\alpha_1,\alpha_2), \\[0.3cm]
    [\phi(\alpha_1,\alpha_2)]\prec P_2= 
  \phi  {\displaystyle\frac{1}{ \lambda}} \ \ln[\cosh(\lambda \alpha_2) +
            \frac{1}{2} \omega \lambda^2 \alpha_1^2 e^{\lambda \alpha_2} + 
(\sinh(\lambda \alpha_2) -  \frac{1}{2} \omega \lambda^2 \alpha_1^2 e^{\lambda \alpha_2}) 
 {\rm C}_\omega(\varphi) 
\\[0.15cm]  \hskip5cm
  + \omega \lambda \alpha_1 e^{\lambda \alpha_2} {\rm S}_\omega(\varphi)]
(\alpha_1,\alpha_2).
   \end{array}$$
The dot stands for the argument of 
the function $\phi=\phi(\; \cdot\; )$. Developing in power series of  $t$ and considering the
first order in  the  first expression   and  the multiplication and derivation  operators
associated to  the  basis
$\varphi^q a_1^r a_2^s$ we can write the  action of  the generators of 
$U_\lambda(\mathfrak{iso}_{\omega}(2))$ over  an arbitrary `function'
$f \in F_\lambda(ISO_{\omega}(2))$
only making  the changes
 \begin{equation} 
\varphi \rightarrow \bar{\varphi}, \qquad
    a_i \rightarrow \frac{\partial }{\partial a_i}\equiv {\partial}_{a_i}.
\end{equation}
In this way we obtain
$$
 \begin{array}{ll} 
    f \prec J =     {\partial}_{\varphi}  f,\\[0.3cm]
    f \prec P_1= { \displaystyle \frac{2 \lambda {\partial}_{a_1}
          e^{\lambda {\partial}_{a_2}} {\rm C}_\omega(\bar{\varphi})+
            (\omega \lambda ^2 \frac{\partial^2}{\partial a_1^2}
                  e^{\lambda {\partial}_{a_2}} 
- {2\sinh(\lambda {\partial}_{a_2})})
  {\rm S}_\omega(\bar{\varphi})}{
2\lambda \cosh(\lambda {\partial}_{a_2})+
     \omega \lambda^3 \frac{\partial^2}{\partial a_1^2}
        e^{\lambda {\partial}_{a_2}}+
          (2\lambda \sinh(\lambda {\partial}_{a_2})-
            \omega \lambda^3 \frac{\partial^2}{\partial a_1^2}
         e^{\lambda {\partial}_{a_2}})  {\rm C}_\omega (\bar{\varphi}) 
                + 2\omega \lambda ^2 {\partial}_{a_1}
       e^{\lambda {\partial}_{a_2}} {\rm S}_\omega(\bar{\varphi})} f },
\\[0.3cm]
   f \prec P_2=   \frac{1}{\lambda}
            \ln[\cosh(\lambda {\partial}_{a_2})+
           \frac{1}{2} \omega \lambda^2 \frac{\partial^2}{\partial a_1^2}
              e^{\lambda {\partial}_{a_2}}+ 
                (\sinh(\lambda {\partial}_{a_2})-
                    \frac{1}{2} \omega \lambda^2 \frac{\partial^2}{\partial a_1^2}
           e^{\lambda {\partial}_{a_2}}   {\rm C}_\omega (\bar{\varphi}) +
  \omega \lambda {\partial}_{a_1}
                 e^{\lambda {\partial}_{a_2}} {\rm S}_\omega (\bar{\varphi})] f.  
 \end{array}
$$
From the above result the  action  of the generator $J$ over  the ordering monomials $\varphi^q
a_1^r a_2^s$ is easily evaluated
\begin{equation}
  (\varphi^q a_1^r a_2^s) \prec J=    q  \varphi^{q-1} a_1^r a_2^s,
\end{equation}
but it does not happen the same with $P_1$ and $P_2$, excepting in the nondeformed case
where the action reduces to
\begin{equation}
    f \prec P_1= \left({\rm C}_\omega(\bar{\varphi}) \frac{\partial}{\partial a_1}- 
{\rm S}_\omega(\bar{\varphi}) \frac{\partial}{\partial a_2}\right)f,    \qquad
    f \prec P_2=  \left( \omega {\rm S}_\omega(\bar{\varphi}) \frac{\partial}{\partial a_1}+
                  {\rm C}_\omega(\bar{\varphi}) \frac{\partial}{\partial a_2}\right)f. 
   \end{equation}
So, when $\lambda=0$ one gets:
\begin{equation} \begin{split}
     (\varphi^q a_1^r a_2^s) \prec P_1= &  r  \varphi^q {\rm C}_\omega(\varphi) a_1^{r-1} a_2^s
         -   s  \varphi^q {\rm S}_\omega(\varphi) a_1^{r} a_2^{s-1}, \\[0.3cm]
     (\varphi^q a_1^r a_2^s) \prec P_1= & 
                        \omega r  \varphi^q {\rm S}_\omega(\varphi) a_1^{r-1} a_2^s
                          +   s  \varphi^q {\rm C}_\omega(\varphi) a_1^{r} a_2^{s-1}.
\end{split}\end{equation}

{\bf 2).-} The  description of the  left coregular module 
 $(F_\lambda(ISO_{\omega}(2)), \succ ,U_\lambda(\mathfrak{iso}_{\omega}(2)) )$
is made in an analogous way. Firstly, considering  Theorem~\ref{crbicross}
one obtains
\begin{equation}\begin{split}
  e^{tJ} \succ & [\phi(\alpha_1, \alpha_2)]=   \phi(\, \cdot \, e^{tJ}) \\[0.15cm]
                    & \left(\frac{2 \lambda\alpha_1 e^{\lambda \alpha_2} {\rm C}_\omega(t)
 + (\omega \lambda ^2 \alpha_1^2 e^{\lambda \alpha_2}
        -2 \sinh(\lambda \alpha_2)) {\rm S}_\omega(t)}
{2\lambda\cosh(\lambda \alpha_2) +
                       \omega \lambda^3 \alpha_1^2 e^{\lambda\alpha_2} +
                    (2\lambda\sinh(\lambda \alpha_2) -
           \omega \lambda^3 \alpha_1^2 e^{\lambda \alpha_2})
                                   {\rm C}_\omega(t) +
              2\omega \lambda^2 \alpha_1 e^{\lambda \alpha_2} {\rm S}_\omega(t)}\right. ,
\\[0.15cm]
                     & \left.  \frac{1}{ \lambda}
                \ln[\cosh(\lambda \alpha_2) +
            \frac{1}{2} \omega \lambda^2 \alpha_1^2 e^{\lambda \alpha_2} +
                    (\sinh(\lambda \alpha_2) -
               \frac{1}{2} \omega \lambda^2 \alpha_1^2 e^{\lambda \alpha_2})
                                   {\rm C}_\omega(t) +
       \omega \lambda \alpha_1 e^{\lambda \alpha_2} {\rm S}_\omega(t)]\right), \\[0.3cm]
  P_1 \succ & [\phi(\alpha_1, \alpha_2)]=  \alpha_1  \phi(\alpha_1, \alpha_2), \\[0.3cm]
  P_2 \succ & [\phi(\alpha_1, \alpha_2)]=   \alpha_2  \phi(\alpha_1, \alpha_2).
\end{split} \end{equation}

With the same arguments that in the case of the right regular co-space
 the  action of  the generators can be written in terms of  the 
multiplication and derivation operators
\begin{equation}
 \begin{split}
   & J \succ f= [\frac{\partial}{\partial \varphi}-
           \bar{a}_1 \frac{1- e^{-2 \lambda \frac{\partial }{\partial a_2}}}{2 \lambda}
              +\omega \bar{a}_2 \frac{\partial}{\partial a_1}-
     \frac{ \lambda}{2} \omega \bar{a}_1 \frac{\partial^2}{\partial a_1^2}] f, \\
   & P_i \succ f=  \frac{\partial}{\partial a_i} f,\quad i=1,2,
 \end{split}
\end{equation}
which allows to obtain the  action of  the generators over  the  basis 
$\varphi^q a_1^r a_2^s$
\begin{equation}
\begin{split}
J \succ (\varphi^q a_1^r a_2^s)= \ q \varphi^{q} a_1^{r} a_2^{s-1}& -
                 \frac{\lambda}{2} \omega r (r-1)
       \varphi^q a_1^{r-1} a_2^{s} + \omega r \varphi^{q} a_1^{r-1} a_2^{s+1} \\
   &  - \frac{1}{2 \lambda} \varphi^{q} a_1^{r+1} a_2^{s}
     + \frac{1}{2 \lambda} \varphi^{q} a_1^{r+1} (a_2+2 \lambda)^{s}.
      s \varphi^{q} a_1^{r} a_2^{s-1}, \\[0.3cm]
P_1 \succ (\varphi^q a_1^r a_2^s)=&  r \varphi^q a_1^{r-1} a_2^s, \\[0.3cm]
P_2 \succ (\varphi^q a_1^r a_2^s)=& s \varphi^q a_1^{r} a_2^{s-1}.
\end{split}
\end{equation}

The  subalgebra  ${\cal A}_{\omega, \lambda}=\langle a_1, a_2 \rangle$ of 
$F_\lambda(ISO_{\omega}(2))$, which is not a  Hopf  subalgebra,
is stable under  the previous action.  The explicit  action  of  the
generators of  $U_\lambda(\mathfrak{iso}_{\omega}(2))$ over  a generic element $\psi(a_1,a_2)$
of  ${\cal A}_{\omega, \lambda}$ in terms of operators adapted to 
the  basis $a_1^r a_2^s$ is
   \begin{equation} \label{planowl}
   \begin{split}
   & J \lact \psi(a_1,a_2)= \left(-
     \bar{a}_1 \frac{1- e^{-2 \lambda \frac{\partial }{\partial a_2}}}{2 \lambda}
              +\omega \bar{a}_2 \frac{\partial}{\partial a_1}-
     \frac{ \lambda}{2} \omega \bar{a}_1 \frac{\partial^2}{\partial a_1^2}\right)
                \psi(a_1,a_2), \\[0.3cm]
& P_i \lact \psi(a_1,a_2)=  \frac{\partial}{\partial a_i} \psi(a_1,a_2),
\qquad i=1,2.
 \end{split}
   \end{equation}
The generators  $a_1$ and $a_2$ of the  co-space
$({\cal A}_{\omega, \lambda}, \lact,  U_\lambda(\mathfrak{iso}_{\omega}(2)))$ verify  the 
commutation relation
   \begin{equation}
             [a_1, a_2]= \lambda a_1
   \end{equation}
defining a noncommutative geometry except in the   nondeformed limit 
where the corresponding classical geometries are recovered.

The  action  of the  Casimir (\ref{casimir}) of $U_\lambda (\mathfrak{iso}_{\omega}(2))$ on
$({\cal A}_{\omega, \lambda}, \lact,  U_\lambda(\mathfrak{iso}_{\omega}(2)))$
is given by 
\begin{equation}\label{casimiriso}
   {\rm C}_{\omega, \lambda} \lact f(a_1, a_2) = 
   \left[ \omega  \frac{\partial^2 }{\partial a_1^2} e^{\lambda \frac{\partial}{\partial a_2}}+
       \frac{4}{\lambda^2}  \sinh^2(\frac{\lambda}{2} \frac{\partial }{\partial a_2})\right]
f(a_1, a_2),
 \end{equation}
or explicitly
\begin{equation}
    {\rm C}_{\omega, \lambda} \lact f(a_1, a_2) = 
 \omega \frac{\partial^2 }{\partial a_1^2}f(a_1, a_2+ \lambda) +
       \frac{1}{\lambda^2}[f(a_1, a_2+ \lambda) + f(a_1, a_2 - \lambda) - 2 f(a_1, a_2)].
\end{equation}
This last expression  shows  the   effect of  the  deformation transforming one
of the derivatives on a finite difference operator.

An interesting problem is  the  solution of  the wave equations associated to this
two--parameter family of Casimir operators (in this context see, for instance Ref.~\cite{olmo96}
and references therein).  Considering the group
$T_{\lambda,2}$ included inside of  ${\cal A}_{\omega, \lambda}$ by means of  the 
exponential map,  the  action  of the Casimir over  the elements $(\alpha_1, \alpha_2)\in
T_{\lambda,2}$ is given by 
\begin{equation} 
{\rm C}_{\omega,\lambda} \lact (\alpha_1, \alpha_2) =
  \left[\omega \alpha_1^2 e^{\lambda \alpha_2}+
   \frac{4}{\lambda^2}\sinh^2(\frac{\lambda}{2} \alpha_2)\right] (\alpha_1, \alpha_2),
 \end{equation}
that suggests to interpret to $(\alpha_1,\alpha_2)$ as a ``plane wave''.

\subsect{Representations of  $U_\lambda(\mathfrak{iso}_{\omega}(2))$}
\label{quantumstandardulambdarepresenyations}

The  representation induced by the  character of  $\cal L$
   \begin{equation} 
1 \dashv (P_1^{n_1} P_2^{n_2})=
      \alpha_1^{n_1}  \alpha_2^{n_2}
\end{equation}
(or, in other words, the  representation  induced by $(\alpha_1, \alpha_2) \in
T_{\lambda,2}$) is obtained as follows. Theorem~\ref{repequiv} gives the   following  result
for  the  action of  the generators in  the  induced representation
\begin{equation} \label{repinduiso}
 \begin{split}
    \phi \dashv J  = &\phi', \\[0.3cm]
    \phi \dashv P_1  = &   \phi \
 \frac{2\lambda \alpha_1 e^{\lambda \alpha_2} {\rm C}_\omega(\varphi) +
    (\omega \lambda ^2 \alpha_1^2 e^{\lambda \alpha_2}
         -2\sinh(\lambda \alpha_2)){\rm S}_\omega(\varphi)}
{ 2\lambda \cosh(\lambda \alpha_2 +
         2\omega \lambda^3 \alpha_1^2 e^{\lambda\alpha_2}+
                    (2\lambda \sinh(\lambda \alpha_2) -
 \omega \lambda^3 \alpha_1^2 e^{\lambda \alpha_2})  {\rm C}_\omega(\varphi) +
      2\omega \lambda ^2 \alpha_1 e^{\lambda \alpha_2} {\rm S}_\omega(\varphi)},\\[0.3cm]
    \phi \dashv P_2 = & \phi \ \frac{1}{\lambda}\ \ln[\cosh(\lambda \alpha_2) +
            \frac{1}{2} \omega \lambda^2 \alpha_1^2 e^{\lambda \alpha_2} \\[0.1cm]
    &  \qquad \qquad +    (\sinh(\lambda \alpha_2) -
      \frac{1}{2} \omega \lambda^2 \alpha_1^2 e^{\lambda \alpha_2}) {\rm C}_\omega(\varphi) +
          \omega \lambda \alpha_1 e^{\lambda \alpha_2} {\rm S}_\omega(\varphi)].
   \end{split}\end{equation}
In the   limit $\lambda\to 0$ we recover the more familiar expressions
 \begin{equation} \begin{split}
  \phi \dashv J = & \phi', \\
   \phi \dashv P_1 = & \phi(\alpha_1 {\rm C}_\omega(\varphi)-
\alpha_2 {\rm S}_\omega(\varphi)), \\[0.3cm]
    \phi \dashv P_2 = & \phi(\omega \alpha_1
{\rm S}_\omega(\varphi)+\alpha_2 {\rm C}_\omega(\varphi)).     
  \end{split}\end{equation}

The local representations   obtained inducing
with the  character of  $U(\mathfrak{so}_\omega(2))$ given by 
\begin{equation} 
J^m \vdash 1= c^m,
\end{equation} 
are
\begin{equation}\begin{split}
  e^{tJ} \vdash & (\alpha_1, \alpha_2)=   e^{tc} \\
        & \times \left( \frac{2\lambda \alpha_1 e^{\lambda \alpha_2} {\rm C}_\omega(t) +
       (\omega \lambda ^2 \alpha_1^2 e^{\lambda \alpha_2}
           -2\sinh(\lambda \alpha_2)) 
{\rm S}_\omega(t)}{2\lambda \cosh(\lambda \alpha_2) +
        \omega \lambda^3 \alpha_1^2 e^{\lambda\alpha_2}+
            (2\lambda \sinh(\lambda \alpha_2) -
   \omega \lambda^3 \alpha_1^2 e^{\lambda \alpha_2}) {\rm C}_\omega(t) +
      2\omega \lambda^2 \alpha_1 e^{\lambda \alpha_2} {\rm S}_\omega(t)},\right. \\
  & \left.  \frac{1}{ \lambda} \ln[\cosh(\lambda \alpha_2) +
      \frac{1}{2} \omega \lambda^2 \alpha_1^2 e^{\lambda \alpha_2}
+  (\sinh(\lambda \alpha_2) -
  \frac{1}{2} \omega \lambda^2 \alpha_1^2 e^{\lambda \alpha_2}) {\rm C}_\omega(t) +
   \omega \lambda \alpha_1 e^{\lambda \alpha_2} {\rm S}_\omega(t)]\right), \\& \\
  P_i \vdash & (\alpha_1, \alpha_2)=  \alpha_i (\alpha_1, \alpha_2),\qquad\qquad i=1,2,
\end{split} \end{equation}
which can be written in an equivalent way using an arbitrary element,
$\psi(a_1,a_2)$, of  $U(\mathfrak{t}_{\lambda,2})$, as 
\begin{equation}\begin{split}
   J\vdash \psi(a_1,a_2)=& \left[c-
   \bar{a}_1 \frac{1- e^{-2 \lambda \frac{\partial }{\partial a_2}}}{2 \lambda}
              +\omega \bar{a}_2 \frac{\partial}{\partial a_1}-
  \frac{ \lambda}{2} \omega \bar{a}_1 \frac{\partial^2}{\partial a_1^2}
                \right] \psi(a_1,a_2), \\
     P_i \vdash \psi(a_1,a_2)= & \frac{\partial}{\partial a_i} \psi(a_1,a_2), \qquad\qquad
i=1,2.  
  \end{split} \end{equation}
Notice that if we take $c=0$, which is equivalent to consider the   character
determined by the  counit of  $U(\mathfrak{so}_\omega(2))$, 
 the  action  (\ref{planowl}) of the  co-space
 $({\cal A}_{\omega, \lambda}, \lact,  U_\lambda(\mathfrak{iso}_{\omega}(2)))$ is recovered.

 The  Casimir action  for the local representation coincides with  those given by 
 (\ref{casimiriso}) since only  the generators $P_i$ are presented.

It is worthy to note that also in this case in the limit $\lambda\to 0$ we obtain
 \begin{equation} \begin{split}
 J   \vdash \psi(a_1,a_2)= & \left[c -
   {a}_1 \frac{\partial}{\partial a_2}
              +\omega {a}_2 \frac{\partial}{\partial a_1} \right] \psi(a_1,a_2), \\[0.3cm]
   P_i \vdash \psi(a_1,a_2)= & \frac{\partial}{\partial a_i} \psi(a_1,a_2), \qquad\qquad
i=1,2    \end{split}\end{equation}
that it is in agreement with the results for local representations \cite{olmo92}.

\sect{Concluding remarks}

The induction procedure that we have used here is not, strictly speaking, a generalization of 
Mackey's induction method for Lie groups. The concept of co-space generalizes in an
algebraical way the concept of $G$--space (being $G$ a transformation group),  and is 
 connected with the induced representations in some sense similar to the nondeformed case.
In Ref.~\cite{olmo00} we presented a more general method that, obviously,  can also be used for
bicrossproduct algebras. However, it is necessary the knowledge of pairs of dual basis of
the corresponding Hopf algebra and its dual, but here this requirement is not necessary in the
method used in this work.

We have made use of the vector fields to compute commutators.  For bicrossproduct Hopf
algebras, like  used in this work $H={\cal K}\bicross {\cal L}$, there is  a connection  between
the representations of $H$  induced by characters of $L$ and  one-parameter flows. This
relation allows to associate  to quantum
bicrossproduct groups  dynamical systems. Here we have only sketched this
situation that  will be  analyzed with more detail  in a forthcoming paper. 

The equations associated to the Casimir operators, as in the nondeformed case, will give the
behaviour of the `deformed' quantum  systems. A procedure for this solution can be found in
\cite{olmo96}. Note that $q$--special polynomials and $q$--functions may appear as solutions of
these $q$--Casimir equations. 

As we mention in Section~\ref{preliminaries}  the bicrossproduct structure  gives  a
$*$--structure  for which the representations are, essentially, unitary. The equivalence of the
induced representations  is given by Theorem~\ref{repequiv}, which establishes a correspondence 
among classes of induced representations and orbits of $L$ under the action of
$K$. This result is analogous to the Kirillov orbits method \cite{kirillov}.  The
problem of the irreducibility of the representations is still open. Partial results for
particular cases have been obtained (see Ref.~\cite{olmo98} for the $(1+1)$ $\k$--Galilei
algebra and \cite{bgst98} for the quantum extended $(1+1)$ Galilei algebra).


\section*{Acknowledgments}
This work has been partially supported by DGES of the Ministerio de Educaci\'{o}n y
Cultura de Espa\~na under Project PB98-0360, and the
Junta de  Castilla y Le\'on (Spain).


\end{document}